\documentclass[12pt]{article}
\usepackage{xcolor}
\usepackage{a4wide}
\usepackage{amsfonts}
\usepackage{amsmath}
\newtheorem{thm}{Theorem}[section]

\newcommand{\dE}{\mathbb{E}}
\newcommand{\dP}{\mathbb{P}}

\def\E{{\mathbb E}}

\newcommand{\ind}{\mbox{1}\kern-.25em \mbox{I}}

\def\build#1_#2^#3{\mathrel{\mathop{\kern 0pt#1}\limits_{#2}^{#3}}}

\def\videbox{\mathbin{\vbox{\hrule\hbox{\vrule height1ex \kern.5em
\vrule height1ex}\hrule}}}

\date{}
\title{A complete characterization of a correlated Bernoulli process}

\author{M. Gonz\'alez-Navarrete\thanks{Departamento de Matem\'atica y Estad\'istica, Universidad de La Frontera, Temuco, Chile. \textbf{manuel.gonzaleznavarrete@ufrontera.cl}}, \  R. Lambert\thanks{Faculdade de Matem\'atica, Universidade Federal de Uberl\^andia, Uberl\^andia, Brasil. \break \textbf{rodrigolambert@yahoo.com.br}}      \ \ and V.H. V\'azquez Guevara\thanks{Facultad de Ciencias F\'isico Matem\'aticas, Benem\'erita Universidad Aut\'onoma de Puebla, Puebla, \break Mexico. \textbf{victor.vazquezg@correo.buap.mx}}}

\begin{document}

\maketitle




\abstract{We present a complete characterization of the asymptotic behaviour of a correlated Bernoulli sequence { which depends on the parameter $\theta \in [0,1]$. A martingale theory based approach will allow} us to prove versions of the law of large numbers, quadratic strong law, law of iterated logarithm, almost sure central limit theorem and functional central limit theorem, in the case $\theta \le 1/2$. For $\theta > 1/2$, we will obtain a strong convergence to a non-degenerated random variable, including a central limit theorem and a law of iterated logarithm for the fluctuations.}


\newcommand{\ABS}[1]{\left(#1\right)} 




\section{Introduction}
The asymptotics for the success rate of a Bernoulli sequence is a pretty well known subject matter in probability and statistics. Its intuitive character and great applicability made it quite popular in the scientific community, and is often used for modelling problems on different areas. 

Since the independent case has already been well studied, many works { have considered} some dependence structure on the source. A particularly interesting case {arises} when the probability of success depends on the number (or rate) of previous successes. In this sense, {we refer to the earliest works \cite{rutherford,woodbury}, in which the authors { deal with the probability distribution of random variables with this kind of structure.}

By following the approach of looking to the previous successes rate to determine the next step's probability distribution, the present paper considers a generalization of the binomial distribution proposed in \cite{DF}. In this case, for $n\geq 0$,the associated random variable $S_n$ counts the number of successes in a correlated Bernoulli sequence, denoted by $\{X_n, n\ge 1\}$ and with conditional probabilities given by
\begin{equation}\label{model}
\dP ( X_{n+1} = 1 | S_n ) = (1-\theta ) p + \theta \frac{S_n}{n},
\end{equation}
where $S_n=\displaystyle\sum_{i=1}^n X_i$ at time $n\geq 0$, and $\dP(X_1=1) = \alpha$. In this sense, we may imagine that $S_n$ represents the number of infected individuals in a population and each person is infected with a probability that depends on previous contagious; or $\{X_n, n\ge 1\}$ represent binary opinions and each individual can be influenced by previous opinions.} { In addition, from a technical point of view, successes probability at time $n\geq 0$ are a weighing of the previous rate of successes and the fixed probability $p$.}

In the case $\alpha=p$, it is possible (as we can see in \cite{heyde2004}) to show that $\E(X_n)= p$ and $\E(S_n)= np$, for any $n\geq 1$, and the limiting behaviour of $S_n$ depends on the parameter $\theta$, in the sense that we could have under or over dispersion, as follows:
\begin{equation}
\label{variance}
\E(S_n - np)^2 \sim \left\{
\begin{array}{rl}
\displaystyle\frac{p(1-p)n}{1-2\theta}     &, \  \mbox{if } \  \theta < 1/2, \\[0.3cm]
p(1-p)n \log (n)     &, \  \mbox{if } \  \theta = 1/2, \\[0.3cm]
\displaystyle\frac{p(1-p)n^{2\theta}}{(2\theta-1)\Gamma(\theta)}     &, \  \mbox{if } \  \theta > 1/2.
\end{array}
\right.
\end{equation}

{ At this same work, author} also proves { versions of the} central limit theorem in the regions $\theta \le 1/2$. Later, {versions of the law of iterated logarithm were demonstrated in \cite{JJQ}}. Moreover, {in\cite{heyde2004}; for $\theta > 1/2$,  was also showed }an almost sure convergence to a random variable, say $L$, for which {it} was possible to obtain the first and second moments. In fact, the papers \cite{Dre,heyde2004} discussed some evidence that $L$ must be not {normally} distributed. Also {at} this region, for finite values of $n$, the random variable $S_n$ could be bimodal, as showed in \cite{Dre}.

In this paper we include a detailed asymptotic {analysis} of the random variable $S_n$, by including functional central limit theorem, almost sure central limit theorem, law of iterated logarithm and convergence of moments in the regions $\theta \le 1/2$. In addition, for $\theta > 1/2$, we obtain a strong convergence to a non-degenerated random variable, showing a dependence of its moments on the initial probability $\alpha$. Moreover, we include a central limit theorem and a law of iterated logarithm for the fluctuations.

The rest of the paper is organized as follows. The next section states the main results. We finish the paper in Section \ref{proofs}, which is dedicated to the proofs.

\section{Main Results}

This section presents the main results of this work. Namely we state a collection of limiting theorems for the number of successes $(S_n)$ on the correlated Bernoulli process conducted by \eqref{model}.


{Our first result deals with the whole range of parameter $\theta$, and shows the almost sure convergence of the successes rate to the fixed probability $p$.}

\begin{thm}
\label{LLN}
For all $\theta\in [0,1)$,
\begin{equation}
\label{LLNequ}
\displaystyle\lim_{n \rightarrow \infty} \frac{S_n}{n}=p \ \text{ a.s}
\end{equation}

\end{thm}



{At the $\theta < 1/2$ setting, we provide the following limit theorems:}

\begin{thm}
\label{T_conv}
If $0 \leq \theta < 1/2$, then, as $n \to \infty$
\begin{itemize}
    \item [i)] We have the distributional convergence in $D([0,\infty[)$ {the Skorokhod space of right-continuous functions with left-hand limits},
\begin{equation}
\label{FCLT-DR}
\left( \sqrt{n}\Big(\frac{S_{\lfloor nt \rfloor}}{\lfloor nt \rfloor}-p\Big), t \geq 0\right) \Longrightarrow \big( W_t, t \geq 0 \big)
\end{equation}
where $\Longrightarrow$ stands for convergence in distribution, and $\big( W_t, t \geq 0 \big)$ is a real-valued centered Gaussian process starting at the origin with covariance given by $\dE[W_s W_t]= \frac{p(1-p)}{(1-2\theta)t} \Big(\frac{t}{s}\Bigr)^{\theta}$, for all $0<s \leq t$.
In particular, $\sqrt{n}\Big(\frac{S_n}{n}-p\Big) \Longrightarrow  N\Big(0,\frac{p(1-p)}{1-2\theta}\Big)$.
 \item[ii)] We have the following almost sure convergence of empirical measures
\begin{equation}
\label{asclt}
\begin{array}{lll}
  \displaystyle\frac{1}{\log n}\sum_{k=1}^n\frac{1}{k}\delta_{\left\{\sqrt{k}\left(\frac{S_k}{k} - p \right) \leq x\right\}}\Rightarrow G   & \text{ a.s} \ ,
\end{array}
\end{equation}
where $\delta_{x}(A)$ equals 1, if $x \in A$, and 0 otherwise. Moreover, $G$ is the Gaussian measure $N\left(0,\frac{p(1-p)}{1-2\theta} \right)$.

  \item[iii)] We obtain the law of iterated logarithm

    \begin{equation}
    \label{LIL1}
        \limsup_{n \to \infty} \pm \frac{S_n-np}{\sqrt{2 \frac{p(1-p)}{1-2\theta} n \log \log \left(n\right)}} = 1 \ \text{ a.s }
    \end{equation}

\item[iv)] It holds the following almost sure convergence
\begin{equation}
\label{moments1}
\begin{array}{lll}
  \displaystyle\frac{1}{\log n}\sum_{k=1}^n k^{r-1}\left(\frac{S_k}{k} - p \right)^{2r} \rightarrow \frac{(p(1-p))^r(2r)!}{2^rr!(1-2\theta)^r }
\end{array}
\end{equation}
   
\end{itemize}
\end{thm}

{We continue by displaying the corresponding asymptotic analysis at the framework $\theta=1/2$}

\begin{thm}
\label{criticalregion}

If $\theta=1/2$, then, as $n \to \infty$, we have
\begin{enumerate}

\item [i)] the distributional convergence in $D([0,\infty[)$,
\begin{equation}
\label{FCLT-CR}
\left( \sqrt{\frac{n^t}{\log n}}\Big(\frac{S_{\lfloor n^t \rfloor}}{\lfloor n^t \rfloor}-p\Big), t \geq 0\right) \Longrightarrow \left(p(1-p) B_t, t \geq 0 \right)
\end{equation}
where $\big( B_t, t \geq 0 \big)$ is a standard Brownian motion.
In particular, \break $\sqrt{\frac{n}{\log n}}\left( \frac{S_n}{n}-p\right) \Longrightarrow  N\left( 0,p(1-p)\right).$

\item [ii)] the almost sure central limit theorem
\begin{equation}
\label{asclt2}
\begin{array}{lll}
  \dfrac{1}{\log\log n}\displaystyle\sum_{k=2}^n\frac{1}{k\log k}\delta_{\left\{\sqrt{\frac{k}{log k}}\left(\frac{S_k}{k} -p\right)  \leq x\right\}}\Rightarrow G  & \text{ a.s} \ ,
\end{array}
\end{equation}
where $G$ is the Gaussian measure $N(0, p(1-p))$.

\item[iii)] the law of iterated logarithm
\begin{equation}
    \label{LIL2}
        \limsup_{n \to \infty} \pm \frac{S_n-np}{\sqrt{2 p(1-p) \cdot  n  \cdot \log n \cdot \log \log ( \log n)}} = 1 \ \text{ a.s }
    \end{equation}

\item[iv)] the almost sure convergence
\begin{equation}
\label{moments2}
\begin{array}{lll}
  \displaystyle\frac{1}{\log \log n}\sum_{k=2}^n \left(\frac{1}{k \log k}\right)^{r+1} k^{r-1} \left(\frac{S_k}{k} -p\right)^{2r}\rightarrow \frac{(p(1-p))^r(2r)!}{2^r r!}
\end{array}
\end{equation}    
    
\end{enumerate}
  \end{thm}

Finally, at the last setting we have the following limit results. Whose deductions are technically more complex

\begin{thm}
\label{T-ASP-SR}
If $\theta> 1/2$, as $n \to \infty$

\begin{enumerate}

\item [i)] we have the almost sure convergence,
\begin{equation}
\label{FSLLN-SR}
\left( n^{1-\theta}\Big(\frac{S_{\lfloor nt \rfloor}}{\lfloor nt \rfloor}-p\Big), t > 0\right) \longrightarrow \Big( \frac{1}{t^{1-\theta}}L, t > 0 \Big)
\end{equation}
where $L$ is a non-degenerated random variable such that
\begin{equation}\label{meanL}
\mathbb{E}[L]=\frac{\alpha- p}{\Gamma(\theta+1)} \  \text{ and } \
\mathbb{E}[L^2]=\frac{\alpha+ (\alpha-p)(1 - 4p)+p \left( \frac{1-2\theta p}{2\theta -1}\right)}{\Gamma(2\theta+1)},
\end{equation}
\item [ii)] the Gaussian fluctuations hold
\begin{equation}
    \label{T_cnv_w}
    \sqrt{n^{2\theta-1}}\left( n^{1-\theta}\Big(\frac{S_n}{n}-p\Big)-L \right) 
\Longrightarrow N\Big(0,\frac{p(1-p)}{2\theta-1}\Big)
 \text{ as }  \ n \to \infty
\end{equation}
\item [iii)] we obtain the law of iterated logarithms for fluctuations
    \begin{equation}
    \label{LIL3}
        \limsup_{n \to \infty} \pm \frac{\sqrt{n^{2\theta-1}}\left( n^{1-\theta}\Big(\frac{S_n}{n}-p\Big)-L \right)}{\sqrt{ \log \log n}} = \sqrt{\frac{2p(1-p)}{2\theta-1}} \ \text{ a.s }
    \end{equation}
    \end{enumerate}

\end{thm}

\section{Proofs}\label{proofs}
The proofs are based on the papers \cite{BV2021,GH,GLV2024}. First of all, we construct a martingale associated to $S_n$. {Thus}, let $\left(\mathcal{F}_n\right)$ {be} the increasing sequence of $\sigma$-algebras $\mathcal{F}_n:=\sigma\left(X_1,\ldots,X_n\right)$, {this definition leads us to} 
\begin{equation}\label{expmov}
\mathbb{E}\left[X_{n+1}|\mathcal{F}_n\right]=\theta\frac{S_n}{n}+\omega \hspace{.5cm}\text{a.s.}
\end{equation}
where $\omega:=(1-\theta)p$. Hence
\begin{equation}\label{expos}
\mathbb{E}\left[S_{n+1}|\mathcal{F}_n\right]=\gamma_n S_n+\omega \hspace{.5cm}\text{a.s.,}
\end{equation}%
where $\gamma_n=1+\frac{\theta}{n}.$
Which leads us to base the asymptotic analysis of the RW on the sequence $(M_n)$, given by $M_0=0$ and for $n\geq 1$  by
\begin{equation} \label{martingale}
M_n=a_n S_n-\omega A_n,
\end{equation}
where the sequence $(a_n)$ is given by $a_1=1$ and for $n\geq 2$ as
\begin{equation}\label{an}
a_n=\prod_{k=1}^{n-1}\gamma_k^{-1}=\frac{\Gamma(n)\Gamma(\theta+1)}{\Gamma(n+\theta)}\sim \frac{\Gamma(1+\theta) }{n^\theta},
\end{equation}%
where $\Gamma$ stands for the Euler gamma function. Moreover, the sequence $(A_n)$ is given by $A_0=0$ and for $n\geq 1$ as $A_n=\sum_{k=1}^n a_k$. Additionally, we observe from \eqref{expos} and {\eqref{martingale}} that almost surely
\begin{eqnarray*}
\mathbb{E}[M_{n+1}|\mathcal{F}_n]=a_{n+1}(\gamma_n S_n+\omega)-\omega A_{n+1}=a_n S_n-\omega A_n=M_n.
\end{eqnarray*}
Thus, $(M_n)$ is a discrete time martingale with respect to the filtration $\left(\mathcal{F}_n\right)$.
From the definition of the proposed martingale given in \eqref{martingale}, we observe that
\begin{eqnarray}
\Delta M_{n} = M_{n}-M_{n-1}=a_{n}\left(S_{n}-\gamma_{n-1} S_{n-1}-\omega\right) = a_{n} \xi_{n}, 
\end{eqnarray}
where; for $n\geq 1$, $\xi_n:=S_n-\mathbb{E}[S_n|\mathcal{F}_{n-1}]=S_{n}-(\omega+\gamma_{n-1} S_{n-1})$ then
\begin{equation} \label{martingale2}
M_n=\sum_{k=1}^n a_k \xi_k.
\end{equation}
Furthermore, note that $\mathbb{E}[X_{n+1}^k|\mathcal{F}_{n}] = \mathbb{E}[X_{n+1}|\mathcal{F}_{n}]$, for all $k \ge 1$.
From this, we find that, a.s.
\begin{equation*}
\mathbb{E}[S_{n+1}^2|\mathcal{F}_{n}]=S_n^2\left(1+\frac{2\theta}{n}\right)+S_n\left(2\omega+\frac{\theta}{n}\right)+\omega.
\end{equation*}

In addition, equation \eqref{expos}, leads us to see that $\mathbb{E}\left[\xi_{n+1}|\mathcal{F}_n\right]=0\hspace{.2cm}\text{a.s.}$, which implies
\begin{equation} \label{segunda}
\mathbb{E}\left[\xi^2_{n+1}|\mathcal{F}_n\right]=\left(\theta\frac{S_n}{n}+\omega\right)\left(1-\left(\theta\frac{S_n}{n}+\omega\right)\right),
\end{equation}
and given that $S_n\leq n$, we find that $\sup_{n\geq 0} \mathbb{E}\left[\xi^2_{n+1}|\mathcal{F}_n\right] < \infty$.
On the same direction, it is possible to show that, almost surely
\begin{equation}\label{cuarta}
\sup_{n\geq 0} \mathbb{E}\left[\xi^4_{n+1}|\mathcal{F}_n\right] < \infty.
\end{equation}


Moreover, the predictable quadratic variation of {$(M_n)$} satisfies, for all $n \geq 1$, that
\begin{equation} \label{procrec}
\langle M \rangle_n = \sum_{k=1}^n \mathbb{E}[\Delta M_k^2| \mathcal{F}_{k-1}] = O(v_n)
 \end{equation}
where $v_n=\sum_{k=1}^n a_k^2.$ Then via standard results on the asymptotic of the gamma function, we conclude that{, as $n$ {goes to infinity, it holds that} }

\begin{enumerate}
\item If $\theta<1/2$ then 
\begin{equation}
\label{vn1}
\frac{v_n}{n^{1-2\theta}}\rightarrow \frac{\Gamma^2(\theta+1)}{1-2\theta}
\end{equation}

\item If $\theta=1/2$ then
\begin{equation}
\label{vn2}
\frac{v_n}{\log n}\rightarrow \frac{\pi}{4}
\end{equation}

\item If $\theta>1/2$ then, from \eqref{an} it is possible to deduce that $(v_n)$ converges into a finite value, more precisely
\begin{equation} \label{vn3}
v_n \rightarrow \sum_{k=0}^ \infty \left(\frac{\Gamma(\theta+1)\Gamma(k+1)}{\Gamma(k+\theta+1)} \right)^2 = \setlength\arraycolsep{1pt}
{}_3 F_2\left(\begin{matrix}1& &1& &1\\&\theta+1&
&\theta+1&\end{matrix};1\right),
\end{equation} 
where the above limit is the generalized hypergeometric function.
\end{enumerate}

\subsection{Proof of Theorem \ref{LLN}}
First, note that
\begin{equation}
    \frac{1}{na_n} = \frac{1}{n} \displaystyle\prod_{k=1}^{n-1}\left(1+\dfrac{\theta}{k}\right)=\displaystyle\prod_{k=1}^{n-1}\left(\dfrac{k+\theta}{k+1}\right)=\displaystyle\prod_{k=1}^{n-1}\left(1-\dfrac{1-\theta}{k+1}\right).
\end{equation}
Moreover {since }$0\leq\theta\leq1$, {we imply that}  $0\leq\dfrac{k+\theta}{k+1}\leq 1$. In addition, $\displaystyle\sum_{k=1}^{\infty}\frac{1-\theta}{k+1}=\infty$ implies that $\lim \frac{1}{na_n}=0$. It follows that $({1}/{na_n})_{n\geq1}$ is non increasing. Let define $N_j=\frac{1}{ja_j}\cdot\Delta M_j$. {Given} that $|\Delta M_n|\leq 1$ we conclude that $(N_j)_{j\geq1}$ is a {  martingale difference sequence} such that $\displaystyle\sum_{j=1}^{\infty}\dE[N_j^2\vert \mathcal{F}_{j-1}]\leq\displaystyle\sum_{j=1}^{\infty}\frac{1^2}{j^2}<\infty$. From Theorem 2.17 in \cite{hall2014martingale}, $\displaystyle\sum_{j=1}^{\infty}N_j$ converges almost surely. { An application of Kronecker's lemma together with the fact that ${n}{a_n}\rightarrow\infty$ let us find that}, $\frac{1}{na_n}\displaystyle\sum_{j=1}^{n}\Delta M_j=\frac{1}{na_n}M_n\overset{a.s}{\rightarrow}0 $. Since $\left|\frac{A_n}{n a_n}-\frac{1}{1-\theta}\right| \sim  \frac{1}{(1-\theta)\Gamma(\theta)}\frac{1}{n^{1-\theta}}$, for $\theta \in [0,1)$, and {by remembering that} $w=(1-\theta)p$, we get that
$$
\lim_{n \to \infty}\frac{1}{na_n}\left(a_nS_n-wA_n\right) = \lim_{n \to \infty}\left(\frac{S_n}{n}-p\right) = 0 \ \ \mbox{a.s.}
$$

\subsection{Proof of Theorem \ref{T_conv}}

To prove (i), from \eqref{segunda} and the definition in \eqref{procrec}, we get that
$$
\langle M \rangle_n  = \E[\xi_{1}^2|\mathcal{F}_0]+\sum_{k=1}^{n-1}a_{k+1}^2\left[w(1-w)+\theta \frac{S_k}{k}(1-2w)-\theta^2\frac{S_k^2}{k^2}\right] 
$$
Now, we use Toeplitz lemma to find that
\begin{eqnarray*}
\lim_{n\to \infty} \ \frac{\langle M\rangle_n}{n^{1-2\theta}} & = & \frac{\Gamma^2(1+\theta)}{1-2\theta}\left(w(1-w)+\theta(1-2w)p-\theta^2p^2\right)\\
& = & p(1-p)\frac{\Gamma^2(1+\theta)}{1-2\theta}
\end{eqnarray*}
Therefore, {by } applying the functional central limit theorem for martingales, given in Theorem 2.5 of \cite{DR}, we obtain that
$$
\lim_{n\to \infty} \ \frac{1}{n^{1-2\theta}}\langle M\rangle_{\lfloor nt\rfloor} = p(1-p)\frac{\Gamma^2(1+\theta)}{1-2\theta}t^{1-2\theta} \ \ \mbox{a.s.}
$$
In order to prove Lindeberg's condition, from \eqref{cuarta} note that for any $\varepsilon > 0$ that
\begin{eqnarray*}
\displaystyle\frac{1}{n^{1-2\theta}}\sum_{k=1}^n\dE [\Delta M_k ^2 \mathbb{I}_{\{|\Delta M_k| > \varepsilon \sqrt{n^{1-2\theta}}\}} \vert \mathcal{F}_{k-1}] &\leq&\displaystyle\frac{1}{n^{2(1-2\theta)}\varepsilon^2}\sum_{k=1}^n\dE [\Delta M_k ^4 \vert \mathcal{F}_{k-1}]\\
\leq \displaystyle\frac{1}{n^{2(1-2\theta)}\varepsilon^2}\sum_{k=1}^n a_k^4\dE [\xi_k ^4 \vert \mathcal{F}_{k-1}]
 &\leq& \displaystyle\frac{16}{n^{2(1-2\theta)}\varepsilon^2}\sum_{k=1}^n a_k^4,
\end{eqnarray*}

Therefore, {since $\frac{n^2 a_n^4}{v_n^2}\rightarrow (1-2\theta)^2$ as $n \rightarrow \infty$ implies that $\frac{1}{n^{1-4\theta}}\sum_{k=1}^n a_k^4$ converges to $\frac{\Gamma(\theta+1)^4}{1-4\theta}$, we obtain that}
\begin{equation*}
\displaystyle\frac{1}{n^{1-2\theta}}\sum_{k=1}^n\dE [\Delta M_k ^2 \mathbb{I}_{\{|\Delta M_k| > \varepsilon \sqrt{n^{1-2\theta}}\}} \vert \mathcal{F}_{k-1}] \rightarrow 0 \text{ as } n \to \infty \text{ in probability}.
\end{equation*}
Then, we conclude that for all $t\ge 0$ and for any $\varepsilon > 0$, 
\begin{equation}
\label{Lindfunctional}
\displaystyle\frac{1}{n^{1-2\theta}}\sum_{k=1}^{\lfloor nt \rfloor}\dE [\Delta M_k ^2 \mathbb{I}_{\{|\Delta M_k| > \varepsilon \sqrt{n^{1-2\theta}}\}} \vert \mathcal{F}_{k-1}] \rightarrow 0,
\end{equation}
as $n \to \infty$ in probability. In addition, note that $\lim_{n\to\infty}\frac{{\lfloor nt \rfloor}a_{\lfloor nt \rfloor}}{n^{1-2\theta}}=t^{1-\theta} \Gamma(\theta+1)$, the definition of \eqref{an} and the fact that 
\begin{equation}
\label{B.1Lem}
\frac{A_n}{n a_n}=\frac{1}{\theta-1}\left(\frac{\Gamma(n+\theta)}{\Gamma(n+1)\Gamma(\theta)}-1\right),
\end{equation}
(by Lemma $B.1$ of \cite{ERWBercu}) imply that
\begin{equation}
\label{Mntfunctional}
 \frac{M_{\lfloor nt \rfloor}}{\sqrt{n^{1-2\theta}}}  = \frac{{\lfloor nt \rfloor}a_{\lfloor nt \rfloor}}{\sqrt{n^{1-2\theta}}} \left( \frac{S_{\lfloor nt \rfloor}}{{\lfloor nt \rfloor}} -p \right) + \frac{p\theta}{\sqrt{n^{1-2\theta}}} \hspace{1cm}\text{a.s.,}
\end{equation}
we conclude via Theorem 2.5 of \cite{DR} that $\left( \sqrt{n} \left( \frac{S_{\lfloor nt \rfloor}}{{\lfloor nt \rfloor}} -p\right), t \geq 0\right) \Longrightarrow \big( W_t, t \geq 0 \big),$ where $W_t = B_t /(t^{1-\theta} \Gamma(\theta+1))$, which completes the proof.

\medskip

Lets prove $(ii)$ by using Lemma 4.1 in \cite{GLV2024}, thus
\begin{eqnarray*}
\sum_{k=1}^{\infty}\frac{1}{v_k}\E\left[|\Delta M_k|^2\mathbb{I}_{\{|\Delta M_K|\geq\varepsilon\sqrt{v_k}\}}|\mathcal{F}_{k-1}\right]  \leq  \frac{1}{\varepsilon^2} \sum_{k=1}^{\infty}\frac{1}{v_k^2}\E\left[|\Delta M_k|^4|\mathcal{F}_{k-1}\right] \\
 \leq  \sup_{k\geq 1}\E[\xi_k^4| \mathcal{F}_{k-1}]\frac{1}{\varepsilon^2}\sum_{k=1}^{\infty}\frac{a_k^4}{v_k^2}  \leq  \frac{C}{\varepsilon^2}\sum_{k=1}^{\infty}\frac{a_k^4}{v_k^2} \sim \frac{C}{\varepsilon^2}\sum_{k=1}^{\infty}\frac{(1-2\theta)^2}{k^2} < \infty,
\end{eqnarray*}
for some constant $C$. The second condition of Lemma 4.1 in \cite{GLV2024} is analogously proved by using $a=2$. Therefore we get that
$$
\frac{1}{\log v_n}\sum_{k=1}^n\left(\frac{v_k-v_{k-1}}{v_k}\right)\delta_{M_k/\sqrt{v_{k-1}}} \Longrightarrow G^* \ \mbox{a.s.}
$$
Since the explosion coefficient is given by $
f_k = \frac{v_k-v_{k-1}}{v_k}=\frac{a_k^2}{v_k} \sim \frac{1-2\theta}{k}$,
{and by observing that} $\log v_n \sim (1-2\theta)\log n$ {we obtain that} $\frac{M_k}{\sqrt{v_{k-1}}} \sim \sqrt{\frac{1-2\theta}{k}}\left(S_k-\frac{kw}{1-\theta}\right)$, {which leads us to} complete the proof.

\medskip

To prove $(iii)$, we first remark that Theorem \ref{LLN} together with \eqref{segunda} implies the following almost sure convergence 
\begin{equation}
\label{varianzaxi2}
\lim_{n\rightarrow \infty} \mathbb{E}\left[\xi^2_{n+1}|\mathcal{F}_n\right] = p(1-p),
\end{equation}
{which jointly with 
$$\sum_{k=1}^{\infty}\frac{a_k^4}{v_n^2}=\frac{[\pi(1-2\theta)]^2}{6} \ ´
$$
and }the law of iterated logarithm for martingales (see \cite{stout} for instance), we get that
$$
\limsup_{n \to \infty}\frac{M_n}{\sqrt{2\log\log n}} = -\liminf_{n \to \infty}\frac{M_n}{\sqrt{2\log\log n}} = \sqrt{p(1-p)} \ ,
$$
which implies that
$$
\limsup_{n \to \infty}\left(\frac{n}{2\log\log n}\right)^{1/2}\left(\frac{S_n}{n}-\frac{wA_n}{na_n}\right) = \sqrt{\frac{p(1-p)}{1-2\theta}} \ .
$$
Therefore, { \eqref{B.1Lem}, conduces us to}
$$
\limsup_{n \to \infty}\left(\frac{n}{2\log\log n}\right)^{1/2}\left(\frac{S_n}{n}-p+\frac{\theta p}{na_n}\right) = \sqrt{\frac{p(1-p)}{1-2\theta}} \ .
$$
We remark that, since $
\displaystyle\lim_{n\to \infty}\left(\frac{n}{2\log\log n}\right)^{1/2}\frac{1}{n^{1+\theta}} =0$, we get that the last term vanishes as $n$ diverges, which completes the proof.

\medskip

Finally, we turn the attention to $(iv)$. The proof is based on Lemma 4.2 from \cite{GLV2024}. Notice that, for all $m\in\mathbb{N}$, we get that $\E[X_{n+1}^m|\mathcal{F}_n]\leq 1$. Moreover
\begin{gather*}
\E[\xi_{n+1}^m|\mathcal{F}_n] =  \E[(X_{n+1}-\E(X_{n+1}|\mathcal{F}_n))^m|\mathcal{F}_n] \\
= \E\left[\sum_{i=1}^m{m\choose i}X_{n+1}^{m-i}\left(\E(X_{n+1}|\mathcal{F}_n)^i(-1)^i\right)\right] \leq  (m+1)\max_{i=0,\cdots , m}{m\choose i} \ ,
\end{gather*}
{that} implies that $\sup_{n\geq 0}\E[|\xi_{n+1}|^n|\mathcal{F}_n] < \infty$.
In addition, we get that $f_n \to 0$, as {$n\rightarrow \infty$}, and $
\frac{M_k}{\sqrt{v_{k-1}}} \sim \sqrt{(1-2\theta)k}\left(\frac{S_k}{k}-p\right)$, which completes the proof.

\subsection{Proof of Theorem \ref{criticalregion}}

We will proceed in a similar fashion as in Theorem \ref{T_conv}. First, due to \eqref{an} and \eqref{vn2} we obtain that 
\begin{equation}\label{ll2}
\frac{a_k^4}{v_k^2}\sim \left(\frac{1}{n \log n}\right)^2,
\end{equation}
which implies that 
\begin{equation} \label{ll3}
\sum_{k=1}^\infty \frac{a_k^4}{v_k^2}<\infty.
\end{equation}

Now, {in order to demonstrate} (i) we note that Lindeberg condition holds from \eqref{ll3}. Then, from the functional central limit theorem for martingales \cite{DR}, the definition of $(M_n)$, convergence \eqref{vn2}, and relation equation \eqref{ll2} we complete the proof.

\medskip

In the sequel, we focus in the proof of (ii). Note that the conditions of Lemma 4.1 in \cite{GLV2024} follows from \eqref{ll3}. In addition, it may be found from the definition of $M_n$, the fact that $\left|\frac{A_n}{n a_n}-2\right|\sim \frac{2}{\sqrt{n\pi}}$, \eqref{vn2} and \eqref{ll2} that $\frac{M_n}{\sqrt{v_{n-1}}} \sim \sqrt{\frac{n}{\log n}}\left( \frac{S_n}{n}-p \right),$ which leads to $\displaystyle\frac{1}{\log \log n}\sum_{k=1}^n \frac{1}{k \log k} \delta_{\sqrt{\frac{k}{\log k}}\left(\frac{S_k}{k} - p \right)} \Rightarrow G  \ \text{ a.s.,}$ where $G$ stands for the $N(0,p(1-p))$ distribution.

\medskip

To prove (iii), we use similar arguments than previous theorem, based on \eqref{ll3}.

\medskip

For (iv), we notice that condition of Lemma 4.2 in \cite{GLV2024} holds in the same manner than in the {$\theta<1/2$ regime}. {From} \eqref{varianzaxi2} and \eqref{ll2}, we obtain that $f_n$ converges to zero as $n \rightarrow \infty$. Hence, we may conclude \eqref{moments2} from the definition of $(M_n)$.
%
%

\subsection{ Proof of Theorem \ref{T-ASP-SR}.}

{For} proving i), note that $\mathbb{E}[\xi^2_{n+1}|\mathcal{F}_n] \leq \theta\frac{S_n}{n}+w \leq \theta+w$ and {$(v_n)$} is a non increasing sequence, then we have that
$$\sup_{n\geq 1} \mathbb{E}[M_{n}^2]\leq (\theta + w) \cdot \setlength\arraycolsep{1pt}
{}_3 F_2\left(\begin{matrix}1& &1& &1\\&\theta+1&
&\theta+1&\end{matrix};1\right) <\infty \ .$$
That is to say, martingale $(M_n)$ is bounded in $L^2$. Thus, it converges in $L^2$ and almost surely to the random variable $M=\sum_{k=1}^\infty a_k \xi_k.$ Now, from \eqref{B.1Lem} we have that
$$
\lim_{n \to \infty} na_n\left( \frac{S_n}{n}-p\right) = M-p\theta \ \ \mbox{a.s.} \Longrightarrow \lim_{n \to \infty}n^{1-\theta}\left(\frac{S_n}{n}-p\right) = L := \frac{M - p\theta}{\Gamma(1+\theta)} \ \mbox{a.s.} 
$$
Moreover, given that $(M_n)$ converges to M in $L^2$ we have that \eqref{B.1Lem} and definition of the limit random variable $L$ guide us to
$$\lim_{n\rightarrow \infty} \mathbb{E}\left[\left(n^{1-\theta}\left(\frac{S_n}{n}-p\right) -L \right)^2 \right]=0,$$
We will find now the first two moments of the limiting random variable $L$ given in \eqref{FSLLN-SR}. For this, let us note that $\mathbb{E}[X_1]=\alpha$. In addition, from \eqref{expos} we have; for $n=1,2,\ldots,$ that
\begin{equation}
\label{lEsn}
    \dE[S_n]=\frac{1}{ a_n}\left(\alpha+\omega\cdot  \displaystyle\sum_{l=1}^{n-1}a_{l+1}\right) = \frac{1}{ a_n}\left(\alpha+\omega\cdot (A_n-1)\right)
\end{equation}
{that }directly {lets us see} that $ \mathbb{E}[M_n]=\alpha-\omega,$ which implies that $\lim_{n\rightarrow \infty} \mathbb{E}[M_n]=\mathbb{E}[M]=\alpha-\omega,$ that leads us to
$$\mathbb{E}(L)=\frac{1}{\Gamma(\theta+1)}\left(\mathbb{E} (M)-p\theta\right)=\frac{1}{\Gamma(\theta+1)}(\alpha-\omega-p\theta) = \frac{\alpha-p}{\Gamma(\theta+1)}$$

We will find the second moment of the limiting random variable $L$. For this, note that
\begin{eqnarray}
\mathbb{E}[M_n^2]&=&a_n^2\mathbb{E}[S_n^2]-2\omega a_n A_n\mathbb{E}[S_n]+\omega ^2A_n^2 \notag\\[0.1cm]
&=&a_n^2\mathbb{E}[S_n^2]-2\omega  A_n(\alpha-\omega)-\omega ^2A_n^2. \label{mnsquaree}
\end{eqnarray}
Moreover $\mathbb{E}[S_{n+1}^2]=g_n\mathbb{E}[S_{n}^2]+h_n$, where, $g_n:=1+\frac{2\theta}{n}$, and $h_n:=\left( 2\omega+\frac{\theta}{n}\right) \mathbb{E}[S_n]+\omega$, for $n \geq 1$. Hence, it may be found recursively; for $n \geq 1$, that
\begin{equation} 
\label{MSN4}
\mathbb{E}[S_{n}^2]=\frac{\Gamma(n+2\theta)}{\Gamma(n)}\Big(\frac{\alpha}{\Gamma(2\theta+1)}+\sum_{k=1}^{n-1}h_k\frac{\Gamma(k+1)}{\Gamma(k+1+2\theta)}\Big).
\end{equation}
However, from \eqref{B.1Lem} and \eqref{lEsn} we may see that
\begin{eqnarray*}
h_n = p(1+2n\omega)+(\alpha-p)\frac{(2n\omega +\theta)}{n a_n}.
\end{eqnarray*}
In addition, \eqref{MSN4} and a repeated application of Lemma $B.1$ of \cite{ERWBercu}, that is, for $b\neq a+1$,
$$\sum_{k=1}^{n-1}\frac{\Gamma(k+a)}{\Gamma(k+b)} = \frac{\Gamma(a+1)}{(b-a-1)\Gamma(b)} \left( 1- \frac{\Gamma(n+a)\Gamma(b)}{\Gamma(n+b-1)\Gamma(a+1)} \right),
$$
{lets} us to observe that
\begin{eqnarray}
\hspace{-1cm}
& & \dE[S_n^2]  =\alpha \frac{\Gamma(n+2\theta)}{\Gamma(n)\Gamma(2\theta+1)}
\notag\\
&+& \frac{p}{\Gamma(n)} \left( \frac{(1-2\omega)}{2\theta-1} \left(\frac{\Gamma(n+2\theta)}{\Gamma(1+2\theta)} -  \Gamma(n+1) \right)  - p \left(2\frac{\Gamma(n+2\theta)}{\Gamma(1+2\theta)} -  \Gamma(n+2) \right) \right) \notag\\
&+& \frac{(\alpha-p)}{\Gamma(n)} \left( (1-2\omega) \left(\frac{\Gamma(n+2\theta)}{\Gamma(1+2\theta)} -  \frac{\Gamma(n+\theta)}{\Gamma(1+\theta)}\right)  - 2p \left((\theta+1)\frac{\Gamma(n+2\theta)}{\Gamma(1+2\theta)} -  \frac{\Gamma(n+\theta+1)}{\Gamma(1+\theta)} \right) \right)
\notag
\end{eqnarray}
%
%
Hence, \eqref{an} and \eqref{mnsquaree} lead us to
\begin{equation}
\label{MM2}
\lim_{n\rightarrow \infty} \dE[M_n^2]= \frac{\Gamma^2(\theta+1)}{\Gamma(2\theta +1)}  \left(  \alpha  +  (\alpha-p)(1 - 4p) + p \left( \frac{1-2\theta p}{2\theta -1} \right) \right) +  \theta p(2(\alpha-p)+ \theta p),
\end{equation}
which finally implies that
\begin{equation}
\dE[L^2]=\frac{1}{\Gamma^2(\theta+1)}\Big(\dE[M^2]-2p\theta\dE[M]+ p^2 \theta^2\Big)
=\frac{\alpha +  (\alpha-p)(1 - 4p) + p \left( \frac{1-2\theta p}{2\theta -1} \right)}{\Gamma(2\theta+1)} .
\notag
\end{equation}
Now, the proof of items ii) and iii) will be based on Lemma 4.3 from \cite{GLV2024}. In this sense, note that \eqref{varianzaxi2} and the bounded convergence theorem imply that $\displaystyle\sum_{k=1}^{\infty}\dE[(\Delta M_k)^2]\sim p(1-p)\Gamma(\theta+1)^2\displaystyle\sum_{k=1}^{\infty}\frac{1}{k^{2\theta}}$. Then, since $\theta> 1/2$, we have that $\displaystyle\sum_{k=1}^{\infty}\dE[(\Delta M_k)^2] < \infty$.
Now, given \eqref{varianzaxi2}, we obtain that:
\begin{equation*}
\begin{array}{ll}
     \displaystyle\sum_{k=n}^{\infty}\dE[(\Delta M_k)^2\vert\mathcal{F}_{k-1}]&\sim\displaystyle\sum_{k=n}^{\infty} p(1-p) a_k^2
    \sim p(1-p)\Gamma(\theta+1)^2\displaystyle\sum_{k=n}^{\infty}\frac{1}{k^{2\theta}}\\[0.4cm]
    &\sim\displaystyle\frac{p(1-p)\Gamma(\theta+1)^2}{(2\theta-1)n^{2\theta-1}}
    \sim \frac{p(1-p) }{(2\theta-1)} n a_n^2 ~a.s
\end{array}
\end{equation*}
By using the bounded convergence theorem, it follows that:
\begin{equation}
\label{r2conv}
    r_n^2 := \sum_{k=n}^{\infty}\dE[(\Delta M_k)^2]\sim \frac{p(1-p) }{(2\theta-1)}n a_n^2~a.s
\end{equation}
Then, conditions a) and a') of Lemma 4.3 in \cite{GLV2024} are satisfied. Hence, have that
\begin{equation*}
     r_n^4\sim \frac{\Gamma(\theta+1)^4(p(1-p))^2 }{(2\theta-1)^2n^{4\theta-2}} \text{ and } (\Delta M_n)^4\leq C_1 a_n^4\sim \frac{C_1\Gamma(\theta+1)^4}{n^{4\theta}},
\end{equation*}
{that guides us to} observe that
\begin{equation*}
\begin{array}{c}
         \displaystyle\frac{1}{r_n^2}\sum_{k=n}^{\infty}\dE\left[(\Delta M_{k+1})^2: \vert\Delta M_{k+1}\vert\geq\varepsilon r_k\right]  \leq \displaystyle\frac{1}{\varepsilon^2r_n^4}\displaystyle\sum_{k=n}^{\infty}\dE\left[\vert\Delta M_{k+1}\vert^4\right] \\[0,4cm]
 \leq \displaystyle\frac{C_1}{\varepsilon^2r_n^4}\sum_{k=n}^{\infty} a_k^4
         \leq\frac{C_2}{r_n^4}\displaystyle\sum_{k=n}^{\infty} \frac{1}{k^{4\theta}}\leq \displaystyle C_3 n^{4\theta-2}n^{1-4\theta}
\end{array}
 \end{equation*}
which implies condition b) of Lemma 4.3 in \cite{GLV2024}. Then, by noticing that $M_n-M = a_n \left(S_n - np -n^{\theta}L\right)$, and using \eqref{r2conv}, the convergence \eqref{T_cnv_w} holds.
Additionally, for $\varepsilon>0$ we have that:
\begin{equation*}
         \frac{1}{r_k}\dE\left[\vert\Delta M_{k+1}\vert: \vert\Delta M_{k+1}\vert\geq\varepsilon r_k\right] \leq\frac{1}{r_k} \frac{1}{\varepsilon^3r_k^3}\dE\left[\vert\Delta M_{k+1}\vert^4\right]
         \leq C_4 k^{4\theta-2}k^{4\theta}
         =\frac{C_4}{k^2}
 \end{equation*}
which implies that condition c) of the same Lemma is satisfied.
In addition, given that $\sum_{k=1}^{\infty} \frac{1}{r_k^4}\dE[(\Delta M_k)^4]<\infty$, we obtain condition d).
Finally, let us denote the martingale difference $d_k:=\frac{1}{r_k^2} \left( (\Delta M_k)^2 - \dE[ (\Delta M_k)^2 \vert\mathcal{F}_{k-1}]\right)$, and observe that
\begin{equation*}
\begin{array}{ll}
         \displaystyle\sum_{k=1}^{\infty}\dE\left[d_k^2 \vert\mathcal{F}_{k-1}\right]  = \displaystyle\sum_{k=1}^{\infty}\frac{1}{r_k^4}\left(\dE[ (\Delta M_k)^4 \vert\mathcal{F}_{k-1}]-\dE^2[ (\Delta M_k)^2 \vert\mathcal{F}_{k-1}]\right)\\[0,4cm]
        \leq \displaystyle\sum_{k=1}^{\infty}\frac{1}{r_k^4}\dE[ (\Delta M_k)^4 \vert\mathcal{F}_{k-1}]
         \leq C_5\displaystyle\sum_{k=1}^{\infty} \frac{a_k^4}{r_k^4}\leq C_6\displaystyle\sum_{k=1}^{\infty} \frac{1}{k^2} < + \infty
\end{array}
 \end{equation*}
 Then, as a consequence of Theorem 2.15 from \cite{hall2014martingale}, condition e) of Lemma 4.3 in \cite{GLV2024} is satisfied, and therefore \eqref{LIL3} holds.


\section*{Acknowledgements}
MGN was partially supported by Fondecyt Iniciaci\'on 11200500. RL is
partially supported by FAPEMIG APQ-01341-21 and RED-00133-21 projects.


\end{document}